\documentclass[10pt,twoside]{article}

\usepackage{epsfig,multicol,bbm}
\usepackage{amsfonts,amsmath,amssymb,amsthm,latexsym}
\usepackage[all]{xy}
\usepackage{hyperref}

\date{}

\begin{document}

\centerline{\bf International Journal of Algebra, Vol. x, 200x, no. xx, xxx - xxx}

\centerline{}

\centerline{}

\centerline {\Large{\bf Ladder operators, Fock-spaces, irreducibility and}}

\centerline{}

\centerline{\Large{\bf group gradings for the Relative Parabose Set algebra}}

\centerline{}

\centerline{\bf {Konstantinos Kanakoglou}}

\centerline{}

\centerline{\small{Instituto de F\'{\i}sica y Matem\'{a}ticas (\textsc{Ifm})}}

\centerline{\small{Universidad Michoacana de San Nicol\'{a}s de Hidalgo (\textsc{Umsnh})}}

\centerline{\small{Edificio C-3, Cd. Universitaria, CP 58040, Morelia, Michoac\'{a}n, \textsc{Mexico}}}

\centerline{{\scriptsize and:}}

\centerline{\small{School of Physics, Nuclear and Elementary Particle Physics Department}}

\centerline{\small{Aristotle University of Thessaloniki (\textsc{Auth}), CP 54124, Thessaloniki, \textsc{Greece}}}

\centerline{\small{kanakoglou@hotmail.com, kanakoglou@ifm.umich.mx}}

\centerline{}

\centerline{\bf {Alfredo Herrera-Aguilar}}

\centerline{}

\centerline{\small{Instituto de F\'{\i}sica y Matem\'{a}ticas (\textsc{Ifm})}}

\centerline{\small{Universidad Michoacana de San Nicol\'{a}s de Hidalgo (\textsc{Umsnh})}}

\centerline{\small{Edificio C-3, Cd. Universitaria, CP 58040, Morelia, Michoac\'{a}n, \textsc{Mexico}}}

\centerline{\small{alfredo.herrera.aguilar@gmail.com}}

\newtheorem{Theorem}{\quad Theorem}[section]

\newtheorem{Definition}[Theorem]{\quad Definition}

\newtheorem{Proposition}[Theorem]{\quad Proposition}

\newtheorem{Corollary}[Theorem]{\quad Corollary}

\newtheorem{Lemma}[Theorem]{\quad Lemma}

\newtheorem{Example}[Theorem]{\quad Example}

\begin{abstract}
The Fock-like representations of the Relative Parabose Set (\textsc{Rpbs}) algebra in a single
parabosonic and a single parafermionic degree of freedom are investigated. It is shown that there is an infinite family (parametrized by the
values of a positive integer $p$) of infinite dimensional, non-equivalent, irreducible representations. For each one of them, explicit expressions
are computed for the action of the generators and they are shown to be ladder operators (creation-annihilation operators)
on the specified Fock-spaces. It is proved that each one of these inf. dim. Fock-spaces is irreducible under the action of the whole algebra or
in other words that it is a simple module over the \textsc{Rpbs} algebra. Finally, $(\mathbb{Z}_{2} \times \mathbb{Z}_{2})$-gradings are introduced
for both the algebra $P_{BF}^{(1,1)}$ and the Fock-spaces, the constructed representations are shown to be
$(\mathbb{Z}_{2} \times \mathbb{Z}_{2})$-graded, $P_{BF}^{(1,1)}$-modules and the relation between our present approach and similar works in the
literature is briefly discussed.
\end{abstract}



\section{Introduction}

The \emph{Relative Parabose Set} algebra (\textsc{Rpbs}) has been introduced by Greenberg and Messiah at their seminal paper \cite{GreeMe}, in which the ``free'' parabosonic and parafermionic algebras were also introduced. It has been historically the only -together with a couple of other models introduced in the same paper- attempt for a mixture of (algebraically) interacting parabosonic and parafermionic degrees of freedom. In \cite{KaDaHa}, \cite{Ya1} the \textsc{Rpbs} algebra is studied in detail, from the mathematical viewpoint, and is shown to be isomorphic to the Universal Enveloping algebra (UEA) of a $(\mathbb{Z}_{2} \times \mathbb{Z}_{2})$-graded, $\theta$-colored Lie algebra. Its super-Hopf algebra structure (i.e. braided group structure) is studied and the properties of some subalgebras are investigated. In \cite{KaDaHa} we introduce the notation $P_{BF}$ for the \textsc{Rpbs} algebra in infinite degrees of freedom and the notation $P_{BF}^{(m,n)}$ for \textsc{Rpbs} in $m$-parabosonic and $n$-parafermionic degrees of freedom.

Our central object of study in this paper will be the $P_{BF}^{(1,1)}$ algebra i.e. the Relative Parabose Set algebra \textsc{Rpbs} in a single parabosonic and a single parafermionic degree of freedom. We will describe it in terms of generators and relations: $P_{BF}^{(1,1)}$ is generated -as an associative  algebra- by the four generators $b^{+}$, $b^{-}$ (both corresponding to the parabosonic degree of freedom) and $f^{+}$, $f^{-}$ (corresponding to the parafermionic degree of freedom). The trilinear relations satisfied by these generators are given as follows
\begin{footnotesize}
\begin{equation}  \label{Rpbsingenerandrel-mixed}
\begin{array} {lll}
\big[ \{ b^{+}, b^{+} \}, f^{-} \big] = 0, & \!\! \big[ [ f^{+}, f^{-} ], b^{-} \big] = 0,  & \!\! \big[ \{ b^{+}, b^{+} \}, f^{+} \big] = 0 = \big[ \{ b^{-}, b^{-} \}, f^{+} \big]  \\
\big[ \{ b^{-}, b^{-} \}, f^{-} \big] = 0,  & \!\! \big[ \{ b^{+}, b^{-} \}, f^{-} \big] = 0,  & \!\! \big[ \{ f^{-}, b^{-} \}, b^{-} \big] = 0 = \big[ \{ f^{-}, b^{+} \}, b^{+} \big]  \\
\big[ \{ f^{-}, b^{+} \}, b^{-} \big] = -2 f^{-},  & \!\!  \{ \{ b^{-}, f^{+} \}, f^{-} \} = 2 b^{-},  & \!\! \big[ \{ f^{+}, b^{+} \}, b^{+} \big] = 0 = \big[ \{ f^{+}, b^{-} \}, b^{-} \big] \\
\big[ \{ b^{-}, f^{-} \}, b^{+} \big] =  2 f^{-},  & \!\! \{ \{ f^{-}, b^{-} \}, f^{+} \} =  2 b^{-},  & \!\! \{ \{ b^{-}, f^{-} \}, f^{-} \} = 0 = \{ \{ b^{-}, f^{+} \}, f^{+} \} \\
\big[ \{ b^{-}, b^{+} \}, f^{+} \big] = 0,  & \!\!  \big[ [ f^{-}, f^{+} ], b^{+} \big] = 0,   & \!\!  \{ \{ b^{+}, f^{+} \}, f^{+} \} = 0 = \{ \{ b^{+}, f^{-} \}, f^{-} \} \\
\big[ \{ f^{+}, b^{-} \}, b^{+} \big] = 2 f^{+},  & \!\!  \big[ \{ b^{+}, f^{+} \}, b^{-} \big] = -2 f^{+},   & \!\!  \{ \{ b^{+}, f^{-} \}, f^{+} \} = 2 b^{+} = \{ \{ f^{+}, b^{+} \}, f^{-} \}
\end{array}
\end{equation}
\end{footnotesize}
together with
\begin{footnotesize}
\begin{equation}  \label{Rpbsingenerandrel-pure}
\begin{array}{llll}
\big[ b^{-}, \{ b^{+}, b^{-} \} \big] = 2 b^{-},  & \!\! \big[ b^{+}, \{ b^{+}, b^{+} \} \big]= 0  & \!\!\big[ b^{+}, \{ b^{-}, b^{-} \} \big] = - 4 b^{-}, & \!\! \big[ f^{-}, [ f^{+}, f^{-} ] \big] = 2 f^{-} \\
\big[ b^{-}, \{ b^{-}, b^{-} \} \big]= 0,  & \!\! \big[ b^{-}, \{ b^{+}, b^{+} \} \big] =  4 b^{+}  &  \!\! \big[ b^{+}, \{ b^{-}, b^{+} \} \big] = - 2 b^{+}, & \!\! \big[ f^{+}, [ f^{-}, f^{+} ] \big] = 2 f^{+}
\end{array}
\end{equation}
\end{footnotesize}
One can easily observe that the relations \eqref{Rpbsingenerandrel-pure} involve only the parabosonic and the parafermionic degrees of freedom separately while the ``interaction'' relations \eqref{Rpbsingenerandrel-mixed} mix the parabosonic with the parafermionic degrees of freedom according to the recipe proposed in \cite{GreeMe}. In all the above and in what follows, we use the notation $[x, y]$ (i.e.: the ``\emph{commutator}'') to imply the expression $xy-yx$ and the notation $\{x, y \}$ (i.e.: the ``\emph{anticommutator}'') to imply the expression $xy+yx$, for $x$ and $y$ any elements of the algebra $P_{BF}$.

Technically, the above definition -stated in terms of generators and relations- is equivalent to saying that $P_{BF}^{(1,1)}$ is isomorphic (as an assoc. alg.) to the quotient of the tensor algebra $\mathbb{T}(b^{\pm}, f^{\pm})$ generated by $b^{+}$, $b^{-}$, $f^{+}$, $f^{-}$ via the two sided ideal (of $\mathbb{T}(b^{\pm}, f^{\pm})$) generated (as an ideal) by relations \eqref{Rpbsingenerandrel-mixed}, \eqref{Rpbsingenerandrel-pure}. Consequently we have to do with an associative, infinite dimensional algebra.

As one can easily figure out, not all of the $32 (= 24 + 8)$ relations presented in \eqref{Rpbsingenerandrel-mixed}, \eqref{Rpbsingenerandrel-pure} are algebraically independent. For example it is easy to see that the first and the third relations in the top row of \eqref{Rpbsingenerandrel-pure} are identified. On the other hand relations $\{ \{ f^{+}, b^{+} \}, f^{-} \} = 2 b^{+}$ and $\{ \{ b^{+}, f^{-} \}, f^{+} \} = 2 b^{+}$ can be added to produce $\big[ [ f^{-}, f^{+} ], b^{+} \big] = 0$ etc. However we keep the relations as presented above and we do not proceed to further simplifying them -although it is an easy task- for three distinct reasons: The first, has to do with the fact this form of relations stems directly from the relations of the general case of $P_{BF}$ (see for example the presentation in \cite{KaDaHa}). The second reason, is that this form of writing the relations will be proved to be particularly convenient when we will proceed, in the next sections, in normal-ordering computations in monomials of the generators. The third reason, has to do with the fact that the form \eqref{Rpbsingenerandrel-mixed}, \eqref{Rpbsingenerandrel-pure} of the relations reveals the $\theta$-color, ($\mathbb{Z}_{2} \times \mathbb{Z}_{2}$)-graded Lie algebraic structure (see \cite{KaDaHa, Ya1}) of the Relative Parabose Set algebra. Although we are not going to make direct use of this structure in this article -in fact we are going to use a rather non-compatible form of grading in the last section- we feel it is worth keeping the relations in this form.

In \cite{GreeMe} a general strategy is proposed for constructing representations of the various -either ``free'' or ``mixed''- paraparticle algebras studied there: Conditions are stated, under which a unique, irreducible representation is singled out for each of these algebras. The authors prove the uniqueness and the irreducibility of these representations for the (free) parabosonic and the (free) parafermionic algebras and rather conjecture the corresponding results for the case of the mixed algebras they study (one of these mixed algebras is \textsc{Rpbs}). They also introduce the \emph{Green ansatz}, a device for proceeding with calculations in the above representations. We will call these representations -whose study has begun in \cite{GreeMe}- from now on Fock-like representations, because of their apparent -as we shall see- similarity with the usual boson-fermion Fock space representation. We must underline at this point that, since the time of \cite{GreeMe}, the Fock-like representations of the Relative Parabose Set algebra \textsc{Rpbs} have never been constructed explicitly, not even for the simplest case of the $P_{BF}^{(1,1)}$ algebra (described by rel. \eqref{Rpbsingenerandrel-mixed}, \eqref{Rpbsingenerandrel-pure}), due to the tremendous computational difficulties inserted by the nature and the number of the trilinear relations.

In \cite{Ya2} the authors proceed in studying the Fock-like spaces of $P_{BF}^{(1,1)}$ following exactly the recipe proposed originally by \cite{GreeMe}.

In this article, we are going to compute explicit expressions for the action of the generators of the $P_{BF}^{(1,1)}$ algebra on its Fock spaces, we will show that these generators constitute a kind of generalized ``ladder'' operators (i.e.: generalized creation-annihilation operators) and consequently we will study the subspaces, the irreducibility and the ($\mathbb{Z}_{2} \times \mathbb{Z}_{2}$)-gradings of the Fock spaces. We will show that there is an infinite family (parametrized by the values of a positive integer $p$) of infinite dimensional, non-equivalent, irreducible, $(\mathbb{Z}_{2} \times \mathbb{Z}_{2})$-graded  representations. The structure of the paper will be as follows:

In Section \ref{YJreview}, the construction of the infinite dimensional Fock-spaces, i.e. the carrier spaces for the Fock-like representations of $P_{BF}^{(1,1)}$, carried out in \cite{Ya2} is reviewed. For each positive, integer value of the parameter $p$, we present an infinite set of vectors constituting a basis for the corresponding space, we introduce some terminology and notation and finally we summarize the main results in the form of a theorem.

In Section \ref{2s}, we proceed in a series of lemmas which enable us to compute (after quite lengthy algebraic calculations inside $P_{BF}^{(1,1)}$) explicit expressions for the actions of the $b^{+}$, $b^{-}$ and $f^{+}$, $f^{-}$ generators on the basis vectors. We prove the generators to be ladder operators (creation-annihilation operators) and we then use our results to prove straightforwardly that the Fock-like representations are irreducible or in other words that the inf. dim. Fock-like spaces constitute simple $P_{BF}^{(1,1)}$-modules. In this way we provide a straightforward verification of the corresponding conjecture stated in \cite{GreeMe}.

Finally, in Section \ref{3s} we introduce ($\mathbb{Z}_{2} \times \mathbb{Z}_{2}$)-gradings for both the algebra $P_{BF}^{(1,1)}$ and its (inf. dim.) Fock-spaces. This grading enables us to describe each one of the Fock-like representations (for any fixed value of the positive integer $p$ parametrizing the representations) as an inf. dim., irreducible, $(\mathbb{Z}_{2} \times \mathbb{Z}_{2})$-graded, $P_{BF}^{(1,1)}$-module. However, we comment on the fact that this grading is different (and incompatible in a sense which we explain) with the one studied in the works \cite{KaDaHa, Ya1}.

We conclude the paper with a discussion of our results in Section \ref{4s}.

We remark that all vector spaces, algebras and tensor products in this article will be considered over the field of complex numbers $\mathbb{C}$ and that the prefix ``super'' will always amount to $\mathbb{Z}_{2}$-graded. Finally we use either $P_{BF}^{(1,1)}$ or \textsc{Rpbs} to denote the Relative Parabose Set algebra in a single parabosonic and a single parafermionic degree of freedom.

\section{Fock-spaces structure for $P_{BF}^{(1,1)}$}    \label{YJreview}

Before reviewing the results of \cite{Ya2}, let us recall a conjecture already known from the beginnings of the study of paraparticle algebras (see \cite{GreeMe, OhKa}).

\textbf{Conjecture}[\emph{Greenberg-Messiah, 1965}, see \cite{GreeMe}]\textbf{:}
\emph{If we consider representations of $P_{BF}^{(1,1)}$, satisfying the adjointness conditions $(b^{-})^{\dagger} = b^{+}$ and $(f^{-})^{\dagger} = f^{+}$, on a complex, infinite dimensional, Euclidean \footnote{Euclidean or pre-Hilbert space, in the sense that it is an inner product space, but not necessarily complete (with respect to the inner product).} space possessing a unique vacuum vector $|0 \rangle$ satisfying $b^{-} |0 \rangle = f^{-} |0 \rangle = 0$, then the following conditions ($p$ may be an arbitrary positive integer)
\begin{equation}  \label{singleoutFock}
\begin{array}{ccccc}
b^{-} b^{+} |0 \rangle = f^{-} f^{+} |0 \rangle = p |0 \rangle  & , \qquad &   & \qquad &
b^{-} f^{+} |0 \rangle = f^{-} b^{+} |0 \rangle = 0
\end{array}
\end{equation}
single out an irreducible representation which is unique up to unitary equivalence.}

In other words the above statement, tells us that for any positive integer $p$ there is an irreducible representation of $P_{BF}^{(1,1)}$ uniquely specified (up to unitary equivalence) by $b^{-} |0 \rangle = f^{-} |0 \rangle = 0$ together with relations \eqref{singleoutFock}. It is these representations which we will call Fock-like representations of $P_{BF}^{(1,1)}$ from now on. We emphasize on the fact that each one of these representations is characterized by the positive integer $p$, in other words the value of $p$ is part of the data which uniquely specify the representation.   \\

In \cite{Ya2} the authors investigate the structure of the carrier space of the Fock-like representation of $P_{BF}^{(1,1)}$. Their results may be summarized in the following theorem (notation due to us):
\begin{Theorem}[\textbf{Fock-spaces structure of} $\mathbf{P_{BF}^{(1,1)}}$] \label{FockspstructPBF11}
The carrier spaces of the Fock-like representations of $P_{BF}^{(1,1)}$, uniquely determined (as representations) under the conditions specified in the above conjecture, are:
\begin{equation}  \label{carrierFocksp}
\bigoplus_{n=0}^{p} \bigoplus_{m=0}^{\infty} \mathcal{V}_{m,n}
\end{equation}
where $p$ is an arbitrary (but fixed) positive integer and the subspaces $\mathcal{V}_{m,n}$ are 2-dimensional except for the cases $m = 0$, $n = 0, p$ i.e. except the subspaces $\mathcal{V}_{0,n}$, $\mathcal{V}_{m,0}$, $\mathcal{V}_{m,p}$ which are 1-dimensional. Let us see how the corresponding vectors look like:

$\blacktriangleright$ \underline{If $0 < m$ and $0 < n < p$}, then the subspace $\mathcal{V}_{m,n}$ is spanned by all vectors of the form
\begin{equation}   \label{spanvect}
\Big| \begin{array}{c}
      m_{1}, m_{2}, ..., m_{l}  \\
n_{0},n_{1}, n_{2}, ..., n_{l}
      \end{array} \Big\rangle \equiv (f^{+})^{n_{0}} (b^{+})^{m_{1}} (f^{+})^{n_{1}} (b^{+})^{m_{2}} (f^{+})^{n_{2}} ... (b^{+})^{m_{l}} (f^{+})^{n_{l}} |0 \rangle
\end{equation}
where $m_{1} + m_{2} + ... + m_{l} = m$ , $n_{0} + n_{1} + n_{2} + ... + n_{l} = n$ and $m_{i} \geq 1$ (for $i = 1, 2, ... , l$), $n_{i} \geq 1$ (for $i = 1, 2, ..., l-1$) and $n_{0}, n_{l} \geq 0$.

For any specific combination of values $(m,n)$ the corresponding subspace $\mathcal{V}_{m,n}$ has a basis consisting of the two vectors
\begin{equation} \label{subspbase}
\begin{array}{ccccc}
| m, n, \alpha \rangle \equiv (f^{+})^{n} (b^{+})^{m} | 0 \rangle   &   &  and  &  &
| m, n, \beta \rangle \equiv (f^{+})^{(n-1)} (b^{+})^{(m-1)} R^{+} | 0 \rangle
\end{array}
\end{equation}
where we use the notation $R^{\eta} = \frac{1}{2} \{ b^{\eta}, f^{\eta} \}$ for $\eta = \pm$. In other words we can always express any vector $\Big| \begin{array}{c}
      m_{1}, m_{2}, ..., m_{l}  \\
n_{0},n_{1}, n_{2}, ..., n_{l}
      \end{array} \Big\rangle$
of the form \eqref{spanvect} as a linear combination of vectors of the form \eqref{subspbase}
\begin{equation}
\Big| \begin{array}{c}
      m_{1}, m_{2}, ..., m_{l}  \\
n_{0},n_{1}, n_{2}, ..., n_{l}
      \end{array} \Big\rangle = c_{1} | m, n, \alpha \rangle + c_{2} | m, n, \beta \rangle
\end{equation}

$\blacktriangleright$ \underline{If $m = 0$ or $n = 0, p \ $}, the vectors $|0, n, \beta \rangle$ and $|m, 0, \beta \rangle$ are (by definition) zero and furthermore the vector $|m, p, \beta \rangle$ becomes parallel to $|m, p, \alpha \rangle$, thus:
\begin{equation}
|0, n, \beta \rangle = |m, 0, \beta \rangle  = 0 \ and \ |m, p, \beta \rangle = \frac{1}{p} |m, p, \alpha \rangle
\end{equation}
Consequently, the corresponding subspaces $\mathcal{V}_{0,n}$, $\mathcal{V}_{m,0}$, $\mathcal{V}_{m,p}$ are 1-dimensional and their bases consist of the single vectors $|0, n, \alpha \rangle$, $|m, 0, \alpha \rangle$, $|m, p, \alpha \rangle$ respectively.

$\blacktriangleright$ \underline{If $n \geq p+1$}, all basis vectors of \eqref{subspbase} vanish.
\end{Theorem}

\textbf{Remark:} Note that, according to the notation {\small $\Big| \begin{array}{c}
      m_{1}, m_{2}, ..., m_{l}  \\
n_{0},n_{1}, n_{2}, ..., n_{l}
      \end{array} \Big\rangle$} introduced in relation \eqref{spanvect}, the vectors of relation \eqref{subspbase}
can be written
{\small
\begin{equation} \label{subspbase2}
\begin{array}{ccccc}
| m, n, \alpha \rangle \equiv (f^{+})^{n} (b^{+})^{m} | 0 \rangle =            \Big| \begin{array}{c}
                                                                                      m  \\
                                                                                      n, 0
                                                                                      \end{array} \Big\rangle  \\
   \\
| m, n, \beta \rangle \equiv (f^{+})^{(n-1)} (b^{+})^{(m-1)} R^{+} | 0 \rangle = \frac{1}{2} \Big| \begin{array}{c}
                                                                                                   m  \\
                                                                                                   n-1, 1
                                                                                                   \end{array}\Big\rangle +
\frac{1}{2} \Big| \begin{array}{c}
                  m-1, 1  \\
                  n-1, 1, 0
                  \end{array} \Big\rangle
\end{array}
\end{equation}
}

\section{Action of the generators and irreducibility of the Fock-spaces}     \label{2s}

\paragraph{$\blacksquare$ \textbf{Construction of ladder operators:}}

We will now proceed in computing the action of the generators $b^{+}, b^{-}$ and $f^{+}, f^{-}$ of $P_{BF}^{(1,1)}$ on the basis vectors of the vector space described in Theorem \ref{FockspstructPBF11}. We will show that the generators are indeed ladder operators on the basis specified in Theorem \ref{FockspstructPBF11}.

$\boxed{\mathbf{1}}$ Let us first start with the computation of the action of the $b^{-}$ generator. We start with a couple of intermediate results:
\begin{Lemma}   \label{interm0}
For the $R^{\eta} = \frac{1}{2} \{ b^{\eta}, f^{\eta} \}$  ($\eta = \pm$) elements (defined in Section \ref{YJreview}), of the $P_{BF}^{(1,1)}$ algebra we have:
\begin{equation} \label{0}
(R^{+})^{2} = 0
\end{equation}
\end{Lemma}
\begin{proof}
{\footnotesize
$$
\begin{array}{c}
\{R^{+},R^{+}\}=2(R^{+})^{2}=\frac{1}{4}\big\{\{f^{+},b^{+}\}, \{f^{+},b^{+}\}\big\}=
\frac{1}{4}\Big(\big\{\{f^{+},b^{+}\}, f^{+}b^{+}\big\}+\big\{\{f^{+},b^{+}\}, b^{+}f^{+}\big\}\Big)= \\
=\frac{1}{4}\Big(f^{+}\big\{b^{+},\{f^{+},b^{+}\}\big\}-\big[f^{+},\{f^{+},b^{+}\}\big]b^{+}+b^{+}\big\{f^{+},\{f^{+},b^{+}\}\big\}-\big[b^{+},\{f^{+},b^{+}\}\big]f^{+}\Big) \\
\end{array}
$$
}
but for the last two summands {\footnotesize $\big\{f^{+},\{f^{+},b^{+}\}\big\}=\big[b^{+},\{f^{+},b^{+}\}\big]=0$} (from \eqref{Rpbsingenerandrel-mixed}) thus:
{\footnotesize
$$
\begin{array}{c}
\{R^{+},R^{+}\}=2(R^{+})^{2}=\frac{1}{4}\Big(f^{+}\big\{b^{+},\{f^{+},b^{+}\}\big\}-\big[f^{+},\{f^{+},b^{+}\}\big]b^{+}\Big)= \\
=\frac{1}{2}\Big(f^{+}(b^{+}R^{+}+R^{+}b^{+})-(f^{+}R^{+}-R^{+}f^{+})b^{+}\Big)
\end{array}
$$
}
from which, after applying the third and fourth rel. of \eqref{Rcomrel}, we get
{\footnotesize
$$
\begin{array}{c}
\{R^{+},R^{+}\}=2(R^{+})^{2}=\frac{1}{2}\Big(f^{+}(b^{+}R^{+}+R^{+}b^{+})-(f^{+}R^{+}-R^{+}f^{+})b^{+}\Big)= \\
=\frac{1}{2}\Big(f^{+}b^{+}R^{+}+f^{+}R^{+}b^{+}-f^{+}R^{+}b^{+}+R^{+}f^{+}b^{+}\Big)= \\
=\frac{1}{2}\Big(f^{+}b^{+}R^{+}+f^{+}b^{+}R^{+}-f^{+}b^{+}R^{+}-f^{+}b^{+}R^{+}\Big)=0
\end{array}
$$
}
which finally completes the proof.
\end{proof}
\begin{Lemma}   \label{interm1}
Starting from relations \eqref{Rpbsingenerandrel-mixed} and \eqref{Rpbsingenerandrel-pure} we have the following relations
\begin{equation}   \label{1}
b^{-} (f^{+})^{k} = (-1)^{k+1}(k-1)(f^{+})^{k}b^{-} + (-1)^{k+1}k(f^{+})^{k-1}b^{-}f^{+}
\end{equation}
\begin{equation}   \label{2}
b^{-} (b^{+})^{m} = \left\{%
\begin{array}{ll}
(b^{+})^{m} b^{-} + m (b^{+})^{m-1}, \ \ m: even    \\
(b^{+})^{m-1} b^{-} b^{+} + (m - 1) (b^{+})^{m-1}, \  m:odd
\end{array}%
\right.
\end{equation}
\begin{equation}    \label{3}
f^{+} (b^{+})^{n} = \left\{%
\begin{array}{ll}
(b^{+})^{n}f^{+}, \ n:even   \\
(b^{+})^{n-1}f^{+}b^{+}, \  n:odd
\end{array}%
\right.
\end{equation}
for the integers $k, m, n \geq 0$.
\end{Lemma}
\begin{proof}
For each one of the above relations induction on $k$, $m$, $n$ respectively.
\end{proof}
Also, proceeding inductively we can show (see also \cite{Ya2}) the relation
\begin{equation}   \label{4}
b^{+} (f^{+})^{k} = (-f^{+})^{k} b^{+} + 2 k R^{+} (f^{+})^{k-1}
\end{equation}
for $k \geq 0$. The simplest, non-trivial case of the above is for $k = 1$ in which case \eqref{4} produces $b^{+}f^{+} = -f^{+}b^{+} +2R^{+}$ which will be proved particularly useful in the sequel.
Now we can arrive at the next proposition
\begin{Proposition}  \label{b-act}
Taking into account that $| m, n, \alpha \rangle \equiv (f^{+})^{n} (b^{+})^{m} | 0 \rangle$ and $| m, n, \beta \rangle \equiv (f^{+})^{(n-1)} (b^{+})^{(m-1)} R^{+} | 0 \rangle$ we have the following expressions for the action of the $b^{-}$ generator on the basis vectors
\begin{small}
\begin{equation}   \label{b-1}
\boxed{
b^{-}\!\!\cdot\!| m, n, \alpha \rangle = \left\{%
\begin{array}{ll}
(-1)^{n+2}m | m-1, n, \alpha \rangle + 2(-1)^{n+1}nm | m-1, n, \beta \rangle, \ \underline{m:even}   \\   \\
(-1)^{n+1} \big( 2n-m-(p-1) \big) | m-1, n, \alpha \rangle + \\ \ \ \ \ \ \ \ \ \ \ \ \ \ \ \ \ \ \ \ \ \ \ \ \ \ + 2(-1)^{n+1}n(m-1) | m-1, n, \beta \rangle, \  \underline{m:odd}
\end{array}
\right.
}
\end{equation}
\end{small}
and
\begin{small}
\begin{equation}   \label{b-2}
\boxed{
b^{-}\!\!\cdot\!| m, n, \beta \rangle = \left\{%
\begin{array}{ll}
-(-1)^{n} | m-1, n, \alpha \rangle + (-1)^{n}\big( 2n-m-p \big) | m-1, n, \beta \rangle, \\
\ \ \ \ \ \ \ \ \ \ \ \ \ \ \ \ \ \ \ \ \ \ \ \ \ \ \ \ \ \ \ \ \ \ \ \ \ \ \ \ \ \ \ \ \ \ \ \ \ \ \ \ \ \ \ \ \ \ \ \ \ \ \  \underline{m:even}   \\   \\
-(-1)^{n} | m-1, n, \alpha \rangle - (-1)^{n}(m-1) | m-1, n, \beta \rangle, \  \underline{m:odd}
\end{array}
\right.
}
\end{equation}
\end{small}
for the integers $0 \leq m$, $0 \leq n \leq p$. In the above we have used $`` \cdot "$ to denote the $P_{BF}^{(1,1)}$-action.
\end{Proposition}
\begin{proof}
\underline{For \eqref{b-1}:} We use successively \eqref{1}, \eqref{2}, \eqref{3} and, in the last stage in order to express the resulting vector as a linear combination of the basis vectors \eqref{subspbase}, we make use of $b^{+}f^{+} = -f^{+}b^{+} +2R^{+}$, conditions \eqref{singleoutFock} and of Lemma \ref{interm0} as well.  \\
\underline{For \eqref{b-2}:} We again make successive use of \eqref{1}, \eqref{2}, \eqref{3} and in the last stage we make use of the relations
\begin{equation} \label{Rcomrel}
\begin{array}{ccccc}
\big[ R^{+}, b^{-} \big] = -f^{+}, &   &   \{ R^{+}, f^{-} \} = b^{+}, &  &   \big[ R^{+}, b^{+} \big] = 0 = \{ R^{+}, f^{+} \}
\end{array}
\end{equation}
which are simply a rewriting of the corresponding relations from \eqref{Rpbsingenerandrel-mixed} and \eqref{Rpbsingenerandrel-pure}, together with $b^{+}f^{+} = -f^{+}b^{+} +2R^{+}$, conditions \eqref{singleoutFock} and Lemma \ref{interm0}.
\end{proof}

$\boxed{\mathbf{2}}$ Let us now proceed with the computation of the action of the $f^{-}$ generator. Let us again start with an intermediate Lemma
\begin{Lemma}
Starting from relations \eqref{Rpbsingenerandrel-mixed} and \eqref{Rpbsingenerandrel-pure} we have the following
\begin{equation}   \label{5}
f^{-}(f^{+})^{m} = -(m-1)(f^{+})^{m}f^{-} + m(f^{+})^{m-1}f^{-}f^{+} - m(m-1)(f^{+})^{m-1}
\end{equation}
\begin{equation}    \label{6}
f^{-} (b^{+})^{n} = \left\{%
\begin{array}{ll}
(b^{+})^{n}f^{-}, \ n:even   \\
(b^{+})^{n-1}f^{-}b^{+}, \  n:odd
\end{array}%
\right.
\end{equation}
for the integers $m, n \geq 0$.
\end{Lemma}
\begin{proof}
The proof of both relations can be obtained inductively on $m$ and $n$ respectively (as previously). Notice that \eqref{6} (and \eqref{3} of Lemma \ref{interm1} as well) are direct consequences of the relations
\begin{equation}   \label{7}
\big[ f^{-}, (b^{+})^{2} ] = 0 = \big[ f^{+}, (b^{+})^{2} \big]
\end{equation}
which can be found among \eqref{Rpbsingenerandrel-mixed}.
\end{proof}
Based on the above we will now have the following proposition
\begin{Proposition}    \label{f-act}
Taking into account that $| m, n, \alpha \rangle \equiv (f^{+})^{n} (b^{+})^{m} | 0 \rangle$ and $| m, n, \beta \rangle \equiv (f^{+})^{(n-1)} (b^{+})^{(m-1)} R^{+} | 0 \rangle$ we have the following expressions for the action of the $f^{-}$ generator on the basis vectors
\begin{small}
\begin{equation}   \label{f-1}
\boxed{
f^{-}\!\!\cdot\!| m, n, \alpha \rangle = n(p+1-n) | m, n-1, \alpha \rangle
}
\end{equation}
\end{small}
and
\begin{small}
\begin{equation}   \label{f-2}
\boxed{
f^{-}\!\!\cdot\!| m, n, \beta \rangle =  | m, n-1, \alpha \rangle + (n-1)(p-n)| m, n-1, \beta \rangle
}
\end{equation}
\end{small}
for the integers $0 \leq m$, $0 \leq n \leq p$ (with $`` \cdot "$ we denote the $P_{BF}^{(1,1)}$-action).
\end{Proposition}
\begin{proof}
\underline{For \eqref{f-1}:} After applying \eqref{5}, \eqref{6} and \eqref{3}, in the last stage in order to express the resulting vector as a linear combination of the basis vectors \eqref{subspbase}, we make use of $\big[ [ f^{-}, f^{+} ], b^{+} \big] = 0$ from \eqref{Rpbsingenerandrel-mixed} together with conditions \eqref{singleoutFock}.   \\
\underline{For \eqref{f-2}:} After successively applying \eqref{5}, \eqref{6} and \eqref{3}, we also employ \eqref{Rcomrel}, conditions \eqref{singleoutFock} and $\big[ [ f^{-}, f^{+} ], b^{+} \big] = 0$ as well.
\end{proof}

$\boxed{\mathbf{3}}$ Let us now proceed with the computation of the action of $b^{+}$, $f^{+}$ generators. It is obviously a much easier task than the previous cases. We readily have the following
\begin{Proposition}  \label{b+f+act}
Taking into account that $| m, n, \alpha \rangle \equiv (f^{+})^{n} (b^{+})^{m} | 0 \rangle$ and $| m, n, \beta \rangle \equiv (f^{+})^{(n-1)} (b^{+})^{(m-1)} R^{+} | 0 \rangle$ we have the following expressions for the action of the $b^{+}$, $f^{+}$ generator on the basis vectors
\begin{equation}  \label{b+12}
\begin{array}{c}
\boxed{b^{+}\!\!\cdot\!| m, n, \alpha \rangle = (-1)^{n} | m+1, n, \alpha \rangle + (-1)^{n-1}2n | m+1, n, \beta \rangle}    \\
     \\
\boxed{b^{+}\!\!\cdot\!| m, n, \beta \rangle = (-1)^{n-1} | m+1, n, \beta \rangle}
\end{array}
\end{equation}
{\small
\begin{equation}   \label{f+12}
\begin{array}{c}
\boxed{
f^{+}\!\!\cdot\!\!|m,n,\alpha\rangle\!=\!\left\{%
\begin{array}{l}
|m,n+1,\alpha\rangle, \ \underline{if \ n \leq p-1}   \\
0, \qquad \qquad \qquad \underline{if \ n \geq p}
\end{array}
\right.
}
   \\   \\
\boxed{
f^{+}\!\!\cdot\!\!|m,n,\beta\rangle\!=\!\left\{%
\begin{array}{l}
|m,n+1,\beta\rangle, \ \underline{if \ n \leq p-1}   \\
0, \qquad \qquad \qquad \underline{if \ n \geq p}
\end{array}
\right.
}
\end{array}
\end{equation}
}
for the integers $0 \leq m$, $0 \leq n \leq p$ (with $`` \cdot "$ we denote the $P_{BF}^{(1,1)}$-action).
\end{Proposition}
\begin{proof}
\underline{For \eqref{b+12}:} To prove each one of these two relations, we first employ \eqref{4} and then we also use relations \eqref{Rcomrel} (for the second one of \eqref{b+12}) together with conditions \eqref{singleoutFock}. \\
\underline{For \eqref{f+12}:} These are the easiest case: they are both obvious by the definition of the basis vectors.
\end{proof}
Proposition \ref{b-act}, Proposition \ref{f-act} and Proposition \ref{b+f+act} give explicit expressions for the actions of all the generators of $P_{BF}^{(1,1)}$ on any vector of the Fock spaces $\bigoplus_{n=0}^{p} \bigoplus_{m=0}^{\infty} \mathcal{V}_{m,n}$ defined in Theorem \ref{FockspstructPBF11}. Thus these three propositions fully define the constructed Fock-like representations (for any positive, integer value of $p$) of the Relative Parabose Set algebra (\textsc{Rpbs}) in a single parabosonic and a single parafermionic degree of freedom.

\paragraph{$\blacksquare$ \textbf{Irreducibility of the $P_{BF}^{(1,1)}$ Fock-like representations:}}
We are now going to prove that the Fock-like representations ($\forall \ p \in \mathbb{N}^{*}$) fully determined by Proposition \ref{b-act}, Proposition \ref{f-act} and Proposition \ref{b+f+act} are irreducible representations or equivalently that the Fock spaces $\bigoplus_{n=0}^{p} \bigoplus_{m=0}^{\infty} \mathcal{V}_{m,n}$ are simple $P_{BF}^{(1,1)}$-modules.
\begin{Theorem} \label{Fock-likesum}
The Fock-like representation of the \textsc{Rpbs} algebra (in a single parab. and a single paraf. degree of freedom) is uniquely identified by the conditions \eqref{singleoutFock} stated at the conjecture at the beginning of Section \ref{YJreview}, given that we have chosen some (arbitrary but fixed) value for the positive integer $p$. The carrier space of this representation is the v.s. $\bigoplus_{n=0}^{p} \bigoplus_{m=0}^{\infty} \mathcal{V}_{m,n}$ described in Theorem \ref{FockspstructPBF11}. The explicit expressions for the actions of the generators $b^{+}, b^{-}, f^{+}, f^{-}$ are given in \eqref{b-1}, \eqref{b-2}, \eqref{f-1}, \eqref{f-2}, \eqref{b+12} and \eqref{f+12}. Furthermore, the above vector space has no invariant subspaces under the above defined $P_{BF}^{(1,1)}$-action, thus it is an irreducible representation or equivalently a simple $P_{BF}^{(1,1)}$-module.
\end{Theorem}
\begin{proof}
The fact that the conditions \eqref{singleoutFock}, stated in the conjecture in the beginning of Section \ref{YJreview}, are sufficient for the unique identification of the representation (in other words: \eqref{singleoutFock} carry sufficient information for the full construction of the Fock-like representation) is proved by the fact that for the computation of the explicit forms of the action of the generators $b^{+}, b^{-}, f^{+}, f^{-}$ -given in Proposition \ref{b-act}, Proposition \ref{f-act} and Proposition \ref{b+f+act}- we have only used the trilinear relations \eqref{Rpbsingenerandrel-mixed} and \eqref{Rpbsingenerandrel-pure} of the algebra, together with the relations \eqref{singleoutFock} stated at the conjecture.

The next figure presents all the subspaces $\mathcal{V}_{m,n}$ whose direct sum constitutes the Fock-space $\bigoplus_{n=0}^{p} \bigoplus_{m=0}^{\infty} \mathcal{V}_{m,n}$ (i.e. the carrier space of the Fock-like representation)
{\scriptsize
\begin{equation}   \label{invariantsubspaces0}
\xymatrix{
    \boxed{\mathcal{V}_{0,0}} \ar@<1ex>[r]^{f^{+}} \ar@<1ex>[d]^{b^{+}} & \boxed{\mathcal{V}_{0,1}} \ar@<1ex>[r]^{f^{+}} \ar@<1ex>[l]^{f^{-}} \ar@<1ex>[d]^{b^{+}} & \ldots \ar@<1ex>[r]^{f^{+}} \ar@<1ex>[l]^{f^{-}} &  \boxed{\mathcal{V}_{0,n}} \ar@<1ex>[r]^{f^{+}} \ar@<1ex>[l]^{f^{-}} \ar@<1ex>[d]^{b^{+}} & \ldots \ar@<1ex>[l]^{f^{-}} & \ldots \ar@<1ex>[r]^{f^{+}} & \boxed{\mathcal{V}_{0,p-1}} \ar@<1ex>[r]^{f^{+}} \ar@<1ex>[l]^{f^{-}} \ar@<1ex>[d]^{b^{+}} & \boxed{\mathcal{V}_{0,p}} \ar@<1ex>[l]^{f^{-}} \ar@<1ex>[d]^{b^{+}} \\
    \boxed{\mathcal{V}_{1,0}} \ar@<1ex>[r]^{f^{+}} \ar@<1ex>[d]^{b^{+}} \ar@<1ex>[u]^{b^{-}} & \boxed{\mathcal{V}_{1,1}} \ar@<1ex>[r]^{f^{+}} \ar@<1ex>[l]^{f^{-}} \ar@<1ex>[d]^{b^{+}} \ar@<1ex>[u]^{b^{-}} & \ldots \ar@<1ex>[r]^{f^{+}} \ar@<1ex>[l]^{f^{-}} &  \boxed{\mathcal{V}_{1,n}} \ar@<1ex>[r]^{f^{+}} \ar@<1ex>[l]^{f^{-}} \ar@<1ex>[d]^{b^{+}} \ar@<1ex>[u]^{b^{-}} & \ldots \ar@<1ex>[l]^{f^{-}} &  \ldots \ar@<1ex>[r]^{f^{+}} & \boxed{\mathcal{V}_{1,p-1}} \ar@<1ex>[r]^{f^{+}} \ar@<1ex>[l]^{f^{-}} \ar@<1ex>[d]^{b^{+}} \ar@<1ex>[u]^{b^{-}} & \boxed{\mathcal{V}_{1,p}} \ar@<1ex>[d]^{b^{+}} \ar@<1ex>[u]^{b^{-}} \ar@<1ex>[l]^{f^{-}} \\
    \vdots \ar@<1ex>[d]^{b^{+}} \ar@<1ex>[u]^{b^{-}} & \vdots \ar@<1ex>[d]^{b^{+}} \ar@<1ex>[u]^{b^{-}} & \ldots  & \vdots \ar@<1ex>[d]^{b^{+}} \ar@<1ex>[u]^{b^{-}}  & \ldots \ar@<1ex>[d]^{b^{+}}  & \ldots  & \vdots \ar@<1ex>[u]^{b^{-}}  & \vdots \ar@<1ex>[d]^{b^{+}} \ar@<1ex>[u]^{b^{-}}    \\
    \boxed{\mathcal{V}_{m,0}} \ar@<1ex>[r]^{f^{+}} \ar@<1ex>[d]^{b^{+}} \ar@<1ex>[u]^{b^{-}}  & \boxed{\mathcal{V}_{m,1}} \ar@<1ex>[r]^{f^{+}} \ar@<1ex>[l]^{f^{-}} \ar@<1ex>[d]^{b^{+}} \ar@<1ex>[u]^{b^{-}} & \ldots \ar@<1ex>[r]^{f^{+}} \ar@<1ex>[l]^{f^{-}} & \boxed{\mathcal{V}_{m,n}} \ar@<1ex>[r]^{f^{+}} \ar@<1ex>[l]^{f^{-}} \ar@<1ex>[d]^{b^{+}} \ar@<1ex>[u]^{b^{-}} & \boxed{\mathcal{V}_{m,n+1}} \ar@<1ex>[r]^{f^{+}} \ar@<1ex>[l]^{f^{-}} \ar@<1ex>[d]^{b^{+}} \ar@<1ex>[u]^{b^{-}} & \ldots \ar@<1ex>[l]^{f^{-}}  & \ldots \ar@<1ex>[r]^{f^{+}} & \boxed{\mathcal{V}_{m,p}} \ar@<1ex>[l]^{f^{-}} \ar@<1ex>[d]^{b^{+}} \ar@<1ex>[u]^{b^{-}} \\
    \vdots \ar@<1ex>[u]^{b^{-}} & \vdots \ar@<1ex>[u]^{b^{-}} &  \ldots \ar@<1ex>[r]^{f^{+}}  & \boxed{\mathcal{V}_{m+1,n}} \ar@<1ex>[r]^{f^{+}} \ar@<1ex>[l]^{f^{-}} \ar@<1ex>[d]^{b^{+}} \ar@<1ex>[u]^{b^{-}} & \ldots \ar@<1ex>[l]^{f^{-}} \ar@<1ex>[u]^{b^{-}} &  \ldots  & \vdots  & \vdots \ar@<1ex>[u]^{b^{-}}    \\
    \vdots & \vdots & \ldots & \vdots & \ldots & \ldots  & \vdots & \vdots
}
\end{equation}
}
From relations \eqref{b-1}, \eqref{b-2}, \eqref{f-1}, \eqref{f-2}, \eqref{b+12} and \eqref{f+12} and the above figure we can see that the generators $b^{+}, b^{-}, f^{+}, f^{-}$ of $P_{BF}^{(1,1)}$ are ``\emph{creation-annihilation}" operators in the above defined space. All subspaces constituting the carrier space $\bigoplus_{n=0}^{p} \bigoplus_{m=0}^{\infty} \mathcal{V}_{m,n}$ (represented above) are connected to each other through repeated use (i.e.: monomials or polynomials) of the generators.

This implies that: given any element of the carrier space $\bigoplus_{n=0}^{p} \bigoplus_{m=0}^{\infty} \mathcal{V}_{m,n}$ we can transform it -after a suitable number of applications of the  ``\emph{annihilation}" operators $b^{-}, f^{-}$- to the ground state $| 0 \rangle$ which comprises a basis of the 1-dim. subspace $\mathcal{V}_{0,0}$. On the other hand, starting from $\mathcal{V}_{0,0}$ we can -applying suitable polynomials of the  ``\emph{creation}" operators $b^{+}, f^{+}$- arrive at any given element of  the carrier space.

Thus: the Fock-like representation presented above, is a cyclic module, generated (as a module) by any of its elements. This implies that it is a simple $P_{BF}^{(1,1)}$-module or equivalently an irreducible representation of the $P_{BF}^{(1,1)}$ algebra, which completes the proof.
\end{proof}

At this point, we feel it is quite interesting to remark that inside each one of the $\mathcal{V}_{m,n}$ (at least for those $\mathcal{V}_{m,n}$ whose dim equals $2$ , in other words $m \neq 0$, $n \neq 0, p$) there are operators (i.e. elements of the algebra $P_{BF}^{(1,1)}$) who act as ``ladder operators" (inside $\mathcal{V}_{m,n}$), in other words their action interchanges between linear independent vectors of $\mathcal{V}_{m,n}$. Such an operator is
{\footnotesize
\begin{equation}
T = \frac{p}{2} \big( R^{+}R^{-} + Q^{+}Q^{-} - N_{b} - \frac{p}{2} \big) - 2 \big( N_{b} + \frac{p}{2} \big) \big( N_{f} - \frac{p}{2} \big) N_{s}
\end{equation}
}
where we have used the notation
{\footnotesize
\begin{equation}  \label{defoper}
\begin{array}{cccc}
N_{b} = \frac{1}{2}\{b^{+},b^{-} \}-\frac{p}{2}, &  N_{f}=\frac{1}{2}[f^{+},f^{-}]+\frac{p}{2}, & R^{\eta}=\frac{1}{2}\{b^{\eta},f^{\eta}\}, & Q^{\epsilon}=\frac{1}{2}\{ b^{-\epsilon},f^{\epsilon}\}
\end{array}
\end{equation}
}
for $\epsilon, \eta = \pm$ and {\footnotesize $N_{s} = \frac{1}{p} \big( N_{f}^{2} -(p+1)N_{f} + f^{+}f^{-} + \frac{p}{2} \big)$}.
Such computations, like the action of $T$ which shows it to be a ``ladder operator" inside each one of $\mathcal{V}_{m,n}$, can be found in \cite{Ya2}.

Note that the results of this section and the final statement of the above theorem, fully verify Greenberg´s original conjecture as this is discussed in \cite{GreeMe} and reviewed in the beginning of this section.

\section{The Fock-like representations as ($\mathbb{Z}_{2} \times \mathbb{Z}_{2}$)-graded modules}     \label{3s}

Each Fock-like representation (i.e. $\forall \ p \in \mathbb{N}^{*}$) of the Relative Parabose Set algebra \textsc{Rpbs}, in a single parabosonic and a single parafermionic degree of freedom, described in the previous sections,
was shown to be an infinite dimensional, irreducible, $P_{BF}^{(1,1)}$-module. However, these representations have further interesting and non-trivial mathematical properties.

In this section, we are going to show that these representations are also $(\mathbb{Z}_{2} \times \mathbb{Z}_{2})$-graded representations. We recall at this point the notion of a $\mathbb{G}$-graded module (or: a $\mathbb{G}$-graded representation) where -in principle- $\mathbb{G}$ may be any group: Given a $\mathbb{G}$-graded assoc. algebra $A = \bigoplus_{g \in \mathbb{G}}A_{g}$, a $\mathbb{G}$-graded v.s. $V = \bigoplus_{g \in \mathbb{G}}V_{g}$ and an $A$-action on $V$ then the module ${}_{A}V$ will be called a $\mathbb{G}$-graded module if the $A$-action is ``compatible" with the group operation, i.e. if
\begin{equation} \label{gradedmod}
A_{g} \cdot V_{h} \subseteq V_{gh}
\end{equation}
$\forall$ $g,h \in \mathbb{G}$ where the $A$-action is denoted with ``$\cdot$" and $gh$ is the corresponding group operation in $\mathbb{G}$.

Now we can state the main result of the present section:
\begin{Proposition} \label{gradstructure-vs.}
The Fock-like representations (for any value of the positive integer $p$) described in Theorem \ref{Fock-likesum}, are $(\mathbb{Z}_{2} \times \mathbb{Z}_{2})$-graded modules.
\end{Proposition}
\begin{proof}
Let us start by describing the structure of a $(\mathbb{Z}_{2} \times \mathbb{Z}_{2})$-graded vector space for the carrier space $\bigoplus_{n=0}^{p} \bigoplus_{m=0}^{\infty} \mathcal{V}_{m,n}$: We will consider the following grading on the basis vectors of Theorem \ref{FockspstructPBF11}
\begin{equation}  \label{gradstructure-v.s.}
\verb"deg"|m, n, \alpha \rangle = \verb"deg"|m, n, \beta \rangle = \big( m \ \textsf{mod} \ 2, \ n \ \textsf{mod} \ 2 \big) \in \mathbb{Z}_{2} \times \mathbb{Z}_{2}
\end{equation}
In other words $\verb"deg" \mathcal{V}_{m,n} = (0,0)$ if $m,n$ are both even, $\verb"deg" \mathcal{V}_{m,n} = (1,1)$ if $m,n$ are both odd and $\verb"deg "\mathcal{V}_{m,n} = (0,1)$ \big( or: $(1,0)$ \big) if $m$ is even and $n$ is odd \big( or: $m$ is odd and $n$ is even \big).

Let us now describe a $(\mathbb{Z}_{2} \times \mathbb{Z}_{2})$-graded structure for the $P_{Bf}^{(1,1)}$ algebra. For this purpose we assign the following grading to the generators of the algebra
\begin{equation}   \label{gradstructure-alg.}
\verb"deg" b^{\pm} = (1,0) \ \ \ \ \verb"deg" f^{\pm} = (0,1)
\end{equation}
The $(\mathbb{Z}_{2} \times \mathbb{Z}_{2})$-grading of the generators described in \eqref{gradstructure-alg.}, produces a well-defined $(\mathbb{Z}_{2} \times \mathbb{Z}_{2})$-grading for the whole algebra $P_{BF}^{(1,1)}$. The reason for this is that the trilinear relations of the algebra $P_{BF}^{(1,1)}$ are \emph{homogeneous}: The term homogeneous here amounts to the fact that the lhs and rhs of each one of the relations \eqref{Rpbsingenerandrel-mixed} and \eqref{Rpbsingenerandrel-pure} acquire the same degree under the assignment \eqref{gradstructure-alg.}.

Now it is easy to see, that the $(\mathbb{Z}_{2} \times \mathbb{Z}_{2})$-graded structures described by \eqref{gradstructure-v.s.} and \eqref{gradstructure-alg.} respectively, satisfy \eqref{gradedmod} i.e.
\begin{equation}  \label{gradstructure-mod.}
\begin{array}{cccc}
b^{+} \cdot \mathcal{V}_{m,n} \subseteq \mathcal{V}_{m+1,n}   &   &   &  b^{-} \cdot \mathcal{V}_{m,n} \subseteq \mathcal{V}_{m-1,n}    \\   \\
f^{+} \cdot \mathcal{V}_{m,n} \subseteq \mathcal{V}_{m,n+1}   &   &   &  f^{-} \cdot \mathcal{V}_{m,n} \subseteq \mathcal{V}_{m,n-1}
\end{array}
\end{equation}
Finally, \eqref{gradstructure-mod.} completes the proof.
\end{proof}
We remark, that in all of the above, the $(\mathbb{Z}_{2} \times \mathbb{Z}_{2})$ group is denoted in the additive notation.

Before closing this section, we feel it is worth commenting shortly on a different approach for describing a suitable $(\mathbb{Z}_{2} \times \mathbb{Z}_{2})$-grading for the Relative Parabose Set algebra \textsc{Rpbs} in a single parabosonic and a single parafermionic degree of freedom. The grading of $P_{BF}^{(1,1)}$ is investigated in the works \cite{KaDaHa, Ya1} but from a quite different viewpoint: It is shown that $P_{BF}^{(1,1)}$ is isomorphic to the Universal enveloping algebra of a $\theta$-colored, $(\mathbb{Z}_{2} \times \mathbb{Z}_{2})$-graded Lie algebra. The grading assigned to the generators (in these works) is
\begin{equation}   \label{gradstructure-alg.-alter}
\verb"deg" b^{\pm} = (1,0) \ \ \ \ \verb"deg" f^{\pm} = (1,1)
\end{equation}
One can readily check that the above assignment, also produces a well-defined $(\mathbb{Z}_{2} \times \mathbb{Z}_{2})$-grading for $P_{BF}^{(1,1)}$ or that in other words relations \eqref{Rpbsingenerandrel-mixed} and \eqref{Rpbsingenerandrel-pure} are homogeneous under \eqref{gradstructure-alg.-alter} as well. However, if we consider the $(\mathbb{Z}_{2} \times \mathbb{Z}_{2})$-grading \eqref{gradstructure-alg.-alter} in connection with the $(\mathbb{Z}_{2} \times \mathbb{Z}_{2})$-grading \eqref{gradstructure-v.s.} for the carrier space then we can straightforwardly check that the Fock-like representation is not a $(\mathbb{Z}_{2} \times \mathbb{Z}_{2})$-graded module. On the other hand, the adoption of \eqref{gradstructure-alg.} does not preserve the $\theta$-colored, $(\mathbb{Z}_{2} \times \mathbb{Z}_{2})$-graded Lie algebra structure \cite{KaDaHa, Ya1} for $P_{BF}^{(1,1)}$ but it has the advantage of generating the structure of a $(\mathbb{Z}_{2} \times \mathbb{Z}_{2})$-graded module for the Fock-like representation studied throughout this article.

\section{Discussion} \label{4s}

In this article, we have studied the Fock-like representations (labeled by an arbitrary positive integer $p$) of the Relative Parabose Set algebra (\textsc{Rpbs})  in a single parabosonic and a single parafermionic degree of freedom. We denoted this algebra as $P_{BF}^{(1,1)}$. We have constructed an infinite family (parametrized by the values of a positive integer $p$) of infinite dimensional, non-equivalent, irreducible, $(\mathbb{Z}_{2} \times \mathbb{Z}_{2})$-graded  $P_{BF}^{(1,1)}$-modules.

In Section \ref{2s}, we computed explicit actions for the generators of the algebra on a given basis of the carrier Fock space for an arbitrary (but fixed) value of the positive integer $p$ which parametrises the family of the representations. We also proceeded in a direct proof of the irreducibility of the representations i.e. of the fact that the carrier spaces do not possess any invariant subspaces under the $P_{BF}^{(1,1)}$-action, verifying thus an old conjecture stated in \cite{GreeMe}. The computational results of this section, i.e. the formulae of Proposition \ref{b-act}, Proposition \ref{f-act}, Proposition \ref{b+f+act}, (apart from having been proved in the text) have also been verified with the help of the  \href{http://homepage.cem.itesm.mx/lgomez/quantum/}{Quantum} \cite{GMFD} add-on for Mathematica 7.0, which is an add-on for performing symbolic algebraic computations, including the use of generalized Dirac notation. What we have actually verified via the use of this package, is that all the formula -given by Proposition \ref{b-act}, Proposition \ref{f-act}, Proposition \ref{b+f+act}- determining the action of the generators of $P_{BF}^{(1,1)}$ on the vectors of the carrier spaces $\bigoplus_{n=0}^{p} \bigoplus_{m=0}^{\infty} \mathcal{V}_{m,n}$ are preserving all of the relations of the algebra as these are presented in \eqref{Rpbsingenerandrel-mixed}, \eqref{Rpbsingenerandrel-pure}

In Section \ref{3s}, we assigned $\mathbb{Z}_{2} \times \mathbb{Z}_{2}$ gradings to both the Fock spaces (for any representation of the above family) and the algebra as well and proved that in this way, the Fock-like representations become $(\mathbb{Z}_{2} \times \mathbb{Z}_{2})$-graded modules. We also commented on the differences and the motivations between our present view of the $\mathbb{Z}_{2} \times \mathbb{Z}_{2}$ grading (of the algebra) of \textsc{Rpbs} and previous works \cite{KaDaHa,Ya1} found in the literature.  \\

Before closing this article, we would like to shortly discuss two possible applications of the family of Fock-like representations constructed in this article.

The first one of them, is of pure mathematical interest: it has to do with the possibility to utilize the representations constructed in this article, in conjuction with the Lie superalgebra realizations constructed in \cite{KaDaHa} in order to produce possibly new representations of any Lie superalgebra $L$ initiating from a given $2d$, matrix, $\mathbb{Z}_{2}$-gr. representation of $L$. Work in this direction has already been in progress and we hope we will be able to report further advances in the near future.

The second one, has to do with a potential physical application of the Fock-like representations, in the extension of the study of a well-known model of quantum optics: the Jaynes-Cummings model \cite{JC} is a fully quantized -and yet analytically solvable- model describing (in its initial form) the interaction of a monochromatic electromagnetic field with a two-level atom. Using the Fock spaces built in this article, we might be able to proceed in a generalization of the above model in the study of the interaction of a monochromatic parabosonic field with a $(p+1)$-level system. The Hamiltonian for such a system might be of the form
\begin{equation}
\begin{array}{c}
\mathcal{H} = \mathcal{H}_{b} + \mathcal{H}_{f} + \mathcal{H}_{interact} = \omega_{b}N_{b}+\omega_{f}N_{f}+\lambda(Q^{+}+Q^{-}) = \\ \\
=\frac{\omega_{b}}{2}\{b^{+},b^{-} \}+\frac{\omega_{f}}{2}[f^{+},f^{-}]+\frac{(\omega_{f}-\omega_{b})p}{2}+\frac{\lambda}{2}\big(\{ b^{-},f^{+}\}+\{ b^{+},f^{-}\}\big)
\end{array}
\end{equation}
where $\omega_{b}$ stands for the energy of any paraboson field quanta (this generalizes the photon, represented by the Weyl algebra part of the usual JC-model), $\omega_{f}$ for the energy gap between the subspaces $\mathcal{V}_{m,n}$ and $\mathcal{V}_{m,n+1}$ (this generalizes the two-level atom, represented by the $su(2)$ generators of the usual JC-model) \footnote{actually $\omega_{b}$ and $\omega_{f}$ might be some functions of $m$ or $n$ or both.} and $\lambda$ some suitably chosen coupling constant. The $\mathcal{H}_{b} + \mathcal{H}_{f}$ part of the above Hamiltonian represents the ``field'' and the ``atom'' respectively, while the $\mathcal{H}_{interact} = \lambda(Q^{+}+Q^{-})$ term represents the ``field-atom'' interaction causing transitions from any $\mathcal{V}_{m,n}$ subspace to the subspace $\mathcal{V}_{m-1,n+1} \oplus \mathcal{V}_{m+1,n-1}$ (absorptions and emissions of radiation). The Fock-like representations, the formulas for the action of the generators and the corresponding carrier spaces of this article, provide a full arsenal for performing actual computations in the above conjectured Hamiltonian and for deriving expected and mean values for desired physical quantities. It remains to proceed in a detailed study of the above model and this will be the subject of some next paper.

{\bf {\small Acknowledgements.}} KK would like to thank the whole staff of \textsc{Ifm}, \textsc{Umsnh} for providing a challenging and stimulating atmosphere while preparing this article. His work was supported by the research project \textsc{Conacyt}/No. J60060. The research of AHA was supported by grants \textsc{Cic} 4.16 and \textsc{Conacyt}/No. J60060; he is also grateful to \textsc{Sni}.

{\bf Received: December, 2010}

\end{document}